\magnification=1200
\centerline{\bf Betti numbers of finitely presented groups and very}
\centerline{\bf rapidly growing functions}
\bigskip
\centerline{\bf  Alexander Nabutovsky$^{a,b,1}$ and Shmuel Weinberger$^c$}
\medskip\noindent
{\bf $^a$} Department of Mathematics, 40 St. George st.,
University of Toronto, Toronto, Ontario,
M5S 2E4, CANADA; alex@math.toronto.edu
\par\noindent
{\bf $^b$} Department of Mathematics, McAllister Bldg.,
The Pennsylvania State
University, University Park, PA 16802, USA; nabutov@math.psu.edu
\par\noindent
{\bf $^1$} Corresponding author. Tel.: 1-416-9784793; fax: 1-416-9784107
\par\noindent
{\bf $^c$} Department of Mathematics, 5734 S. University Avenue, University of Chicago,
Chicago, IL 60637-1514, USA; shmuel@math.uchicago.edu.
\bigskip
{\bf Abstract.} Define the length of a finite presentation of a group $G$
as the sum of lengths of all relators plus the number of generators.
How large can be the $k$th Betti number $b_k(G)=$ rank $H_k(G)$ providing
that $G$ has length $\leq N$ and $b_k(G)$ is finite?
We prove that for every $k\geq 3$ the maximum $b_k(N)$ of $k$th Betti numbers
of all such groups is an extremely rapidly growing function of $N$. It grows
faster that all functions previously encountered in Mathematics (outside
of Logic) including non-computable functions (at least those 
that are known to us).
More formally, $b_k$ grows as
the third busy beaver function that measures the maximal productivity
of Turing machines with $\leq N$ states that use the oracle for the halting
problem of Turing machines using the oracle for the halting problem of
usual Turing machines.
\par
We also describe the fastest possible growth of a sequence of finite
Betti numbers of a finitely presented group. In particular, it cannot grow
as fast as the third busy beaver function but can grow faster than the
second busy beaver function that measures the maximal productivity of Turing
machines using an oracle
for the halting problem for usual Turing machines. We describe a natural
problem about Betti numbers of finitely presented groups such that
its answer is
expressed by a function that grows as the fifth busy beaver function.
\par
Also, we outline a construction of a finitely presented
group all of whose homology groups are either ${\bf Z}$ or
trivial such that its Betti numbers form
a random binary sequence. 
\medskip\noindent
{\bf Keywords:} Homology groups of finitely presented groups, Betti numbers,
non-recursive functions, random binary sequences, busy beaver function.
\bigskip\noindent
{\bf 0. Introduction.}
\medskip
In [1] G. Baumslag,E. Dyer and C. Miller gave an almost complete
characterisation of all possible
sequences of homology groups
of finitely presented groups. (Their work was motivated by earlier ideas of
D. Kan and W. Thurston [4].) For example, they have
shown that a sequence $H_1,H_2,H_3,\ldots$
of countably generated torsion-free abelian groups
is a sequence
of all homology groups of a finitely presented group
if and only if 1) $H_1$ and $H_2$ are finitely generated;
and 2) This sequence admits a recursive presentation.
(The notion of recursive presentation
of a sequence of abelian groups can be informally
explained as follows: This is a sequence of countable
presentations
of groups $H_i$ such that there
exists a computer program listing every
relation in all presentations of
the groups $H_i$ (in an arbitrary order). This program
works an infinitely long time, writing from time to time a relation.)
We present a more detailed introduction to results
and methods of [1] in the next section.
\par
Yet this characterization of sequences of homology groups
of finitely presented groups
is not effective enough to make obvious (at least for us)
the answers
for many natural questions about
homology groups of finitely presented groups.
For example, define the $i$th
Betti number of a finitely presented group $G$ as the
rank of the tensor product of $H_i(G)$ with ${\bf R}$.
If this
tensor product is not finitely generated,
we can either define $b_i(G)$ as $\infty$ or regard
$b_i(G)$ as undefined. In the first case we regard
Betti numbers as a function from ${\bf N}$ to ${\bf N}\bigcup \{\infty\}$,
in the second case we can regard Betti numbers
of $G$ as a partial function from ${\bf N}$ to ${\bf N}$, where the
term ``partial" means that the domain of this function is
a subset of ${\bf N}$.
\par
Now we can ask:
how fast such a partial functions can grow?
How fast can it grow if its domain is ${\bf N}$ (that is, all
Betti numbers are finite)?
How large
can be $b_k(G)$  when $k$ is fixed and it is known that
$G$ has a finite presentation of length not exceeding
some (variable) $N$?
\par
In this paper we provide complete answers for these
questions.
These answers are given as Theorems 3.1 and 3.2 in section 3.
To state these results
we need to introduce notions of Turing machines of order $k$
and of $k$th busy beaver functions for every
$k=1,2,\ldots$.
These notions are discussed in the section 2. Here
we would like only to note that the sequence of Betti
numbers of a finitely presented group cannot grow
arbitrarily fast but even in the case when all Betti numbers
of a finitely presented group $G$ are finite, this sequence
can grow more rapidly than any computable function as well
as all known to us non-computable functions previously
encountered in Mathematics (outside of Mathematical Logic).
\par
We deduce these results from Proposition 3.3 and its Corollary 3.3.1 providing
us with a method to effectively realize some sequences
of extended natural numbers as sequences of
Betti numbers of a finitely presented group. 
As another application of Corollary 3.3.1 we indicate how one
can construct an explicit finitely presented group such
that each of its homology groups is either trivial or
isomorphic to ${\bf Z}$ but the sequence of its Betti numbers
is a random sequence of $0$'s and $1$'s (Theorem 4.1).
\par
In the last section
we indicate that a natural question involving
Betti numbers of finitely presented groups has an answer expessed by a
function that grows even more rapidly than functions that appear in
Theorems 3.1, 3.2.
More precisely, this
function grows as the fifth
busy beaver function
(Theorem 5.1).
\par
We are not aware of any natural mathematical
problems that lead to
functions that grow much more rapidly than the fifth busy beaver
function. One
possible source of such problems is ergodic theory (or dynamical
systems),
where one studies outcomes of infinite processes. It is
possible that some natural problems in these areas 
can be stated
only using
predicates with quantifiers with respect to functional
variables, and lead to functions such that
the problem of their computation  belongs to
non-trivial degrees of unsolvability in Kleene's analytic
hierarchy (cf. [7] for an introduction to
Kleene's hierarchies). Yet we do not have any concrete ideas
in this direction.
\medskip\noindent
{\bf 1. Homology groups of finitely presented groups.}
\medskip
Let $G$ be a finitely presented group. It is well-known
that there exists a unique (up to homotopy equivalence) CW
complex denoted $BG$ or $K(G,1)$
such that its fundamental group is isomorphic to
$G$ and all its other homotopy groups vanish. The homology
groups of $BG$ are called {\it homology groups of} $G$.
As usual, the rank of $H_n(G)$ is called {\it the $n$th Betti
number of $G$}.
\par
One possible way to construct $BG$
is the following. First, realize
$G$ as the fundamental group of a finite $2$-complex $K_2$.
(This complex has one $0$-dimensional cell. Its $1$-cells
correspond to generators of $G$, and its $2$-cells correspond
to relators of $G$.) Then one kills all generators of 
$\pi_2(K_2)$ by adding (possibly infinitely many) $3$-cells,
obtainining a $3$-complex $K_3$, and further proceeds
inductively killing on step $i$ all generators of
$\pi_i(K_i)$ by adding $(i+1)$-dimensional cells.
Note that for every $i$ $K_i$ is
a $i$-dimensional CW-complex naturally
included in $K_{i+1}$. One then defines $BG$ as the union
$\bigcup_i K_i$.
\par
This description of $BG$ implies that $H_1(G)$ is just
the abelianization of $G$, $G/[G,G]$ and, thus, is
a finitely generated abelian group. All $2$-cells of $BG$
are already in $K_2$. Therefore $H_2(G)$ is a finitely
generated abelian group. Yet note that 
we added possibly infinitely many $3$-cells, $4$-cells, etc. during our
construction.
Therefore a priori $H_3(G), H_4(G)$, etc. do not need to be
finitely generated. Indeed, J. Stallings [10] constructed
examples of finitely presented groups with infinitely
generated third homology groups. 
\par
A lot of information about homology groups of
finitely presented groups can be found in [1]. There
the authors used the following construction:
Any finitely presented group $G$ can be embedded into an
{\it acyclic} finitely presented group $A_G$. Moreover,
given a finite presentation of $G$ one can explicitly
construct a finite presentation of $A_G$ and the embedding.
(Recall that a group is called acyclic if all its homology
groups are trivial.) This result can be combined with the 
classical theorem of G. Higman: There exists an universal
finitely presented group $U$ such that every countable
recursively presented group $G$ can be effectively
embedded into $U$. Here the effectiveness of embeddability
means that there exists a Turing machine (=an algorithm, a 
computer program) that finds for every generator of $G$
its image under the embedding in $U$. A recursively
presented group is a group with a finite or infinite
countable set of generators and
either a finite or an infinite recursively enumerable 
set of relations. ``Recursive" means here that these
relations are being enumerated by a Turing machine (=by
a computer program): Think about a computer program that
types from time to time a new relation and works infinitely
long. The resulting infinite list of relations will be
an infinite recursively enumerable set of relations.
\par
Embedding the universal Higman group into a finitely
presented acyclic group we obtain a universal acyclic
finitely presented group $A$. Now for every recursively
presented group $A$ we can effectively construct its 
{\it suspension} $SG=A*_GA$. Here we take two copies
of $A$ and embed $G$ into them in the identical way as the
composition of the Higman embedding of $G$ into $U$ and the
embedding of $U$ into $A$.
The term ``effective" means that
there exists an algorithm constructing a recursive
presentation
of the suspension if a recursive presentation of $G$
is given. The output of this algorithm is either a finite
presentation of the suspension if $G$ is finitely generated,
or a Turing machine (=an algorithm) enumerating all
relations of the suspension if $G$ is infinitely generated.
Using the Meyer-Vietoris
exact sequence one immediately sees that for every $i>0$
$H_i(G)=H_{i+1}(A*_GA)$. Note that even if $G$ is infinitely
generated, then its suspension is finitely generated.
If $G$ is finitely generated, then its suspension is
finitely presented. (The idea to use such a group-theoretic suspension
to ``lift" the dimension of a homology group of a group appeared 
already in [4]).
\par
Iterating this construction we obtain the {\it double
suspension of} $G$ $S^2G=A*_{A*_GA}A$.
Note that this group is always finitely presented, and there
exists an algorithm producing a finite presentation
of this group from a given finite presentation of $G$.
For every $i$ $H_{i+2}(S^2G)$
is isomorphic to
$H_i(G)$. In particular, if $G$ is a recursively presented
{\it abelian} group then the third homology group
of the double suspension of $G$ is isomorphic to $G$.
Thus, any recursively presentable
abelian group can be realized as the third homology group of
a finitely presented group. (Vice versa,
the construction of $BG$ outlined above implies that
the third homology group of a finitely presented group
is a recursively presentable abelian group; see [1]
for details.) Further iterating
the suspension construction one can introduce
iterated suspensions $S^kG$ for every $k>2$. All these groups will
be finitely presented, even if $G$ has an infinite
set of generators. Moreover, for
every $i$ $H_{i+k}(S^kG)=H_i(G)$,
and $H_1(G)=G$, if $G$ is abelian.
Thus, in this way one can realize any recusively presented abelian
group as the $k$th homology group of a finitely presented group for any $k>2$.
\par
Further, consider
a recursively
presented sequence of recursively presented abelian groups
with untangled recursive presentations.
(A sequence of recursively presented groups is called
recursively presented if the set of all relations
is a recursively enumerable subset of the set of all words
in all generators of all these groups. In less formal
terms this means that there exists a computer program that
works infinite time, and that writes from time to time
a relation in one of these abelian groups, so that
eventually every relation of every of these groups will
be written down. A finite presentation of an abelian group
is called untangled if for every $l$ first $l$ relations
form a basis of the vector space spanned by these relations.)
The authors of [1] show that if the first two groups 
in this sequence are finitely generated then this sequence
is the sequence
of homology groups $H_1(G),H_2(G),\ldots$ of some finitely
presented group $G$.
Further, it had been shown in [1] that if
a recursively presented abelian group is torsion-free then
one can effectively replace any given recursive presentation
of this group by an untangled recursive presentation.
Thus, all sequences of torsion-free
homology groups of finitely presented groups are
characterised as follows: The class of such sequences coincides
with the class of recursively presented sequences of recursively presented abelian groups, where first two groups
are finitely generated. Moreover, there exists an
algorithm that for every recursively
presented sequence of torsion-free abelian groups constructs
a finitely presented group $G$ such that the groups from the sequence are
isomorphic to $H_3(G), H_4(G), \ldots$.
\par
Note that the same ideas were used in our paper [6] to prove a somewhat
stronger result (Theorem 13.2): If $X$ is any simplicial cell
complex  with computable cell structure then its double suspension
is homotopy equivalent to $K(\pi,1)^+$ for some finitely presented
group $\pi$ that can be explicitly constructed from an algorithm
describing the cell structure of $X$. (Here $^+$ means the Quillen
$^+$-construction with respect to a normal subgroup of $X$.) This
result, which will not be used in the present paper, implies that one
has a wide freedom of realization not only for cohomology groups
of finitely presented groups but also for their
cohomology operations, and more exotic cohomology theories.
\medskip\noindent
{\bf 2. Some facts from computability theory.}
\medskip
In this section we review some well-known facts from recursion theory
that are relevant for our purposes. Our exposition is intended for
readers with very limited previous knowledge of Mathematical Logic.
Most of these facts can be found in [7] or [9], but we hope that a short
self-contained summary is helpful.
\par
For a reader not familiar with Mathematical Logic note
that Turing machines compute exactly the same class of functions
as computer programs in any contemporary programming language such as
C, PASCAL, FORTRAN, etc. providing that this language is stripped of all
data types but integers, and all limitations for the length of numbers,
arrays, etc. (if any) are removed. Such functions are called {\it computable}
or ${\it recursive}$.
A recipee for their computation is called {\it an algorithm}.  A simple
example of a function that is {\it not} computable is Rado's busy
beaver function that is defined as follows: For each $n>1$
consider the (finite) set of all Turing machines
with $\leq n$ states that
eventually stop, when they start the computation with the empty tape.
For each of these machines consider the number of steps of work of the
Turing machine before it stops. Take the maximal of these numbers.  
It can be regarded as a function of $n$. Denote this function $B(n)$.
It is called Rado's busy beaver function. This function dominates every
computable function. (We say that a function $f:{\bf N}\longrightarrow {\bf N}$ dominates
a function $g:{\bf N}\longrightarrow {\bf N}$ if for all sufficiently large $n$
$f(n)>g(n)$.) Of course, one can define a similar function using programs
written in one of programming languages of length $\leq n$ instead of Turing
machines with $\leq n$ states.
\par
 To see that these funcions
dominate every computable function note that
, if $\phi$ is computable, then it can be computed
by a program of constant length (or note that there exists a Turing
machine with $Const$ states, computing $\phi(n)$ for every $n$ written
on tape.) Since $n$ can be represented by a bit sequence of
length $\leq\log n +1$, it is not difficult to see that $\phi(n)$ can be
computed by a program of length $\leq\log n + const$ (or
by a Turing machine with, say $[n/4]$ states (in fact, even much less).
Almost all of these states are required to produce $n$ $1$'s coding $n$ on the
tape.)  Now consider the program that computes $\phi(n)$ in time $T$
and then does $\phi(n)-T$ empty steps before stopping (if $T>\phi(n)$.
Otherwise it stops immediately.) It is clear that for all sufficiently
large $n$ its stopping time is at least $\phi(n)$ but less than $B(n)$.
\par
But why cannot one use the definition of $B(n)$ to design
an algorithm that computes this function? One needs to consider the list
of all Turing machines with $\leq n$ states (or all valid computer
programs in the chosen programming language of length $\leq n$),
find out which of them stop, run all of them, and find the maximal stopping
time. Here all steps with the exception of the second step clearly
can be implemented by means of an algorithm. Thus, we demonstrated the
validity of
the following celebrated Turing theorem: There is no algorithm deciding
whether or not a given Turing machine halts when it starts its computation
with the empty tape.
\par
 (In other words, the halting problem for Turing machines is
algorithmically unsolvable.)
In fact, we have demonstrated more: One can find an algorithm
computing $B(n)$ for any given value of $n$ using an {\it oracle}
deciding whether or not a given Turing machine halts. Vice versa,
if one has an oracle computing $B(n)$ for every given $n$ one can
decide whether or not a given Turing machine halts. Indeed, it is sufficient
to compute $B(n)$, where $n$ is the number of states in the given
Turing machine, and then run it for $B(n)$ steps. If it did not
stop by that time, it will never stop.
\par
Two algorithmic problems are said to have the same (Turing)
degree of unsolvability
if each of them can be solved using an oracle answering the other. 
Turing degrees of unsolvability are equivalence classes of problems
with respect to this equivalence relation.
Thus, computation of the busy beaver function and deciding whether
or not a given Turing machine halt have the same degree of unsolvability
denoted by ${\bf 0'}$. (${\bf 0}$ denotes the degree of unsolvability
of all problems that can be solved by means of an algorithm.) There are
many problems in ``mainstream" mathematics that belong to the
degree of unsolvability ${\bf 0'}$: the word problem and the
triviality problem for finitely presented groups, diffeomorphism problem
for closed smooth manifolds of dimension $\geq 4$, etc., but we will
see soon other problems of even higher degree of unsolvability.
\par
Consider now Turing machines that use an oracle computing $B(n)$ for
any given $n$. Clearly, they are more powerful than usual Turing machines.
For examples, such Turing machines can decide whether or not a  given ``usual"
Turing machine halts. One can define a busy beaver function for such machines
in the same way as
it had been defined for usual Turing machines. More generally,
one can give the following formal definitions:
\medskip\noindent
{\bf Definitions.} A {\it Turing machine of order $1$} is the usual Turing
machine. A {\it Turing machine of order $k$}, $(k>1)$, is a Turing machine
that uses the oracle solving the halting problem for all Turing machines
of order $(k-1)$. The {\it $k$-th busy beaver function} $B_k(n)$ is defined
as the maximal time of work of a Turing machine of order $k$ with $\leq n$
states. (The $n$th state is used only to stop; the machine uses a separate
tape for the oracle; 
we consider only machines that start their work with the
empty tape and eventually halt in forming our maximization.)
In particular, $B_1(n)$ coincides
with Rado's ``busy beaver function" introduced in [7] (see also [3]).
\medskip\noindent
{\bf Theorem 2.1.}
\par\noindent
(a) The halting problem for Turing machines
of order $k$ and the problem of computing $B_k$ belong to the
same degree of unsolvability denoted ${\bf 0^{(k)}}$; these degrees
of unsolvability for different values of $k$ are distinct. 
\par\noindent
(b) For every $k$ and every
$k$-computable function $\phi$  for all sufficiently large $n$
$B_{k+1}(n)>B_k(\phi(n))$.
\medskip
This theorem generalizes the discussion above that corresponds
to the case $k=1$ and can be proven exactly in the same way. (Or see
[9] for a formal proof of (a) and [3] for a formal proof of
(b) in the case $k=1$ that immediately generalizes for all values of $k$.)
\par
It is known that every predicate can be written in a prenex form,
where all quantifiers occur at the beginning of the formula.
Let a predicate be in the prenex form.
Assume that there are $(n-1)$ changes of types of quantifiers
(from the universal to existential or vice versa) in the formula.
Then the predicate is called a $\Sigma_n$-predicate, if the formula
starts from the existential quantifier, and a $\Pi_n$-predicate, if
the formula starts from the universal quantifier. For example,
$\forall x\forall y P(x,y)$ is a $\Pi_1$-predicate, and $\exists x\forall y
\forall z\forall u\exists v\exists w P(x,y,z,u,v,w)$ is a $\Sigma_3$-predicate.
In the present paper
we will consider only first-order arithmetic
predicates where all quantifiers are
applied to variables, and all variables are interpreted as variables with
values in ${\bf N}$. It is known that one can express
the fact that a Turing machine of order $k$ halts when it starts its
work with input $n$ on work tape as a $\Sigma_k$-predicate with one
free (=non-quantified) variable $n$ (cf. [9]). In fact, there is an
algorithm that assigns to a given Turing machine of order $k$ such
a predicate. The first quantified variable in the first block of
existential quantifiers in such a predicate is interpreted as the
halting time of the Turing machine of order $k$. As the result,
the verification of general $\Sigma_k$ (or $\Pi_k$) arithmetic
predicates is as difficult as the halting problem for Turing machines of order $k$. In fact
these two problems belong to the same degree of unsolvability ${\bf 0^{(k)}}$.
\par
A subset of ${\bf N}$ (or ${\bf N}^k$) is called {\it recursively enumerable}
if it is the range of a computable function from ${\bf N}$ to ${\bf N}$ (or ${\bf N}^k$).
An equivalent definition is that a set $A$ is recursively enumerable
if there exists an algorithm (=a Turing machine, a computer program)
that lists all elements of $A$ in some
order. (This algorithm is allowed to work infinite time.) This algorithm
is called {\it an enumeration} of $A$. 
Another equivalent definition is that a set $A$ is recursively enumerable if
it is a domain of a partial computable function; in other
words $A$ is recursively enumerable if
there exists a Turing machine computation of which halts if
and only if the input is in $A$.
A recursively enumerable set can be presented in a finite form by a
Turing machine (or a computer program) enumerating it.
\par
Now one can pose
the following algorithmic problems:
\par\noindent
1. Decide whether or not a given recurively enumerable set is finite;
\par\noindent
2. Decide whether or not a given recursively enumerable set is cofinite.
(That is, decide whether or not the complement of the recursively enumerable
set is finite.)
\par
The first problem is denoted $Fin$ and is known to be in ${\bf 0''}$; the
second problem is denoted $Cof$, and is known to be in ${\bf 0'''}$
(cf. [7], [9]). The upper bounds for complexity of these problems
follow from the possibility to rewrite them as arithmetic predicates
with two blocks of quantifiers (for $Fin$) or three block of quantifiers (for
$Cof$). The proofs of the lower bounds for complexity
 are more delicate. In particular,
the proof of the lower bound for $Cof$ uses the method of moving markers
that will be described in the next section.
\medskip\noindent
{\bf 3. Betti numbers of finitely presented groups}
\medskip
In this section we will prove our main results.
\par
First note that the construction of $BG=K(G,1)$ described in section 1 implies
that for every $n$ we can represent the $n$-skeleton of $K(G,1)$ as the
union of ascending sequence of finite complexes $K_{n,j}$, which can be
effectively constructed. As the result, the Betti numbers of $G$, $b_G(n)$, regarded
as a function of the dimension can be represented as the double limit
$\lim_{j\longrightarrow\infty}\lim_{i\longrightarrow\infty}b_G(i,j,n)$.
Here, $b_G(i,j,n)$ increases with respect to $j$ and decreases with respect to $i$.
The variable $j$ corresponds to $n$-cells of $K(G,1)$ that are being added
all the time and potentially make the $n$th Betti number bigger; the variable
$i$ correspond to the $(n+1)$-dimensional cells that are being added all
the time (to $K_{n,j}$) and potentially make the $n$th Betti  number smaller.
The function $b_G(i,j,n)$, of course, measures the $n$th Betti number
of the intermediate $(n+1)$-dimensional finite cell
complexes that arize as approximations to $K(G,1)$.
It is well-known that a limit of a sequence of computable functions can
be computed by a Turing machine of order $2$, and a double limit can
be computed by a Turing machine of order $3$ (cf. [9]). Thus, we obtain
the first assertion in the following theorem:
\medskip\noindent
{\bf Theorem 3.1.}\par\noindent
(a) The sequence of Betti numbers of any finitely presented group
regarded as a function of the dimension can be computed by a Turing
machine of order $3$.\par\noindent
(b) Let function $f:{\bf N}\longrightarrow
{\bf N}\bigcup\infty$ be any partial function on ${\bf N}$
that can be computed
using a Turing machine of order  $3$ and such that $f(1)$
and $f(2)$ are defined. Then
there exists a finitely presented 
group $G$ such that for every $n$ the $n$th Betti number  $b_n$ of $G$ is finite if and only if
$f(n)$ is defined (i.e. finite), and if $f(n)$ is finite then $b_n\geq f(n)$.
\par
Moreover, there exists an algorithm that constructs a finite presentation
of such a group $G$ starting from a Turing machine of order $3$ computing $f$
as the input data.
\bigskip\noindent
{\bf Corollary 3.1.1.} There exists a finitely presented group $G$ such that
all its Betti numbers are finite, but $b_k(G)>B_2(k)$. On the other hand
for any finitely presented group $G$ and any computable functions
$f_1, f_2:{\bf N}\longrightarrow {\bf N}$ $B_3(n)>f_1(b_{f_2(n)}(G))$
for all sufficiently large $n$.
\medskip\noindent
Indeed, $B_2$ can be computed by a Turing machine of order $3$ (that can use
the oracle solving the halting problem for Turing machines of order $2$).
\medskip
Let $b_n(N)$ denote the maximum of the $n$th Betti numbers among all
finitely presented groups with a finite $n$th Betti number that admit
a finite presentation of length $\leq N$. (The length of a finite
presentation is defined as the sum of lengths of all relators plus the
number of generators in the finite presentation.)
\medskip\noindent
{\bf Theorem 3.2.}
Let $k\geq 3$ be any natural number.
There exist computable functions $f_1, f_2$ such that
for every $n$ $b_k(n)\leq B_3(f_1(n))$ and $b_k(f_2(n))\geq B_3(n)$.
\medskip\noindent
The first of two inequalities in the text of Theorem 3.2 follows from the fact
that $b_k$ of a finitely presented group can be computed by a Turing machine of order $3$. It is easy
to see that the number of states of this machine can be effectively bounded in terms of the
length of a given finite presentation of $G$. The second inequality follows from the second assertion of Theorem 3.1:
The halting time of every Turing machine of order $3$ with $\leq N$ states that halts can be majorized by the $k$th
Betti number of a finitely presented group. Since this group can be effectively constructed, the length of its
finite presentation is effectively bounded in terms of $N$. 
\par
Thus, it remains to prove
the second and the third assertions of Theorem 3.1. 
\medskip\noindent
{\bf Proof of Theorem 3.1:}
To prove Theorem 3.1 recall that
according to [1] every recursively presented sequence of
countably generated torsion-free abelian groups can be represented as the sequence of homology groups of a finitely presented group
providing that the first two abelian groups in this sequence are
finitely generated.
Moreover, if our goal is only to realize these groups as the third, fourth, etc.
homology groups of a finitely presented group $G$, so that the first two homology
groups of $G$ are trivial, then we can find such a $G$ by means of an
algorithm.
Furthermore, it is obvious that one can make the first two Betti numbers of
a finitely presented group arbitrarily large
just by forming the free product of this group with ${\bf Z}^N$ for a
sufficiently large $N$. 
\par
Let $\{I_i\}$ be
a recursive sequence of recursively enumerable sets. (This
means that there exists an algorithm that for each $i$
constructs an enumeration of $I_i$.) Consider
an infinite sequence of recursively presented abelian groups
$A_i$ with abelian generators $x_1, x_2, \ldots$ and
relations $x_j=0$ if an only if $j\in I_i$. Clearly,
this is a recursively presented sequence of abelian groups.
Therefore these groups can be effectively
realized as homology groups
of a finitely presented group. Thus, there exists a
finitely presented group $G$ such that its Betti numbers $b_3, b_4, \ldots$ are cardinalities of the {\it complements}
of sets $I_i$. 
\par
Now it is clear that in order to complete the proof of
Theorem 3.1 it is sufficient to construct an algorithm
that works as follows. For each $n$ this
algorithm provides an enumeration of a 
recursively enumerable set $E_n$ such that its complement is
finite if and only if $f(n)$ is defined. Morover, if $f(n)$
is defined then the cardinality of the complement of 
$E_n$ must be greater than or equal to $f(n)$.
\par
To achieve this goal, first note that 
it is known how to construct
a predicate $P(T)=\exists n_1 \exists n_2 \forall m \exists k Q(n_1, n_2,m,k)$
for each Turing
machine $T$ of order $3$,
so that $T$ halts if and only if $P(T)$ holds (cf. [9]). Here
the meaning of $n_1$ is the
halting time for $T$, $n_2$ codes the computation by $T$ and the oracle
information (=the list of $0$'s and $1$'s coding whether or not first
several Turing machines of order $2$ halt). Further,
$k$ and $m$ code auxilliary
variables required to express that the information obtained from the oracle
is, indeed, what it is supposed to be.
\par
Note that we use here the existence
of an effective  bijection between ${\bf N}$ and a Cartesian product of several
copies of ${\bf N}$. For example, $\phi(n_1, n_2)=(2n_1-1)2^{n_2-1}$ is
a bijection between ${\bf N}\times {\bf N}$ and ${\bf N}$. This bijection can
also be used to
replace the existential quantifiers with respect to $n_1$ and to $n_2$
by one existential quantifier with respect to $N=\phi(n_1,n_2)\geq n_1$). 
We will denote the corresponding
predicate equivalent to $Q(n_1,n_2,m,k)$ by $QQ(N,m,k)$. The minimal
value of $N$ for which $\forall m\exists k QQ(N,m,k)$ is true
is greater than or equal to
the smallest value of $n_1$ for which
$\exists n_2\forall m\exists k Q(n_1, n_2, m, k)$ is true.
\par
Now Theorem 3.1 immediately follows from the next Proposition:
\bigskip\noindent
{\bf Proposition 3.3.} There exists an algorithm that produces for
every given arithmetic
$\Sigma_3$ predicate $P=\exists n\forall m\exists k QQ(n, m, k)$
an enumeration of a recursively enumerable set such
that its complement is finite if and only if $P$
is true, and if the complement is finite then its
cardinality is equal to the minimal value of $n$ for
which $\forall m\exists k\ QQ(n,m,k)$.
\medskip\noindent
Our discussion above (based on results of [1]) and Proposition 3.3
imply the following result:
\medskip\noindent
{\bf Corollary 3.3.1.} (i) There exists an algorithm that
for each $l>2$ and each arithmetic
$\Sigma_3$ predicate $P=\exists n\forall m\exists k\ QQ(n,m,k)$ constructs
a finitely presented group $G$ such that $b_l(G)$
is finite if and only if $P$ is true, and if $P$ is true,
then $b_l(G)$ is equal to the minimal value of $n$ for
which $\forall m\exists k\ QQ(n,m,k)$.
\par\noindent
(ii) For every $\Sigma_3$ predicate
$P=\exists n\forall m\exists k\ QQ(n,m,k,l)$ there exists a finitely presented
group $G$ such that for every $l>2$ $b_l(G)$ is equal to the minimal
$n$ such that \par\noindent
$\forall m\exists k\ QQ(n,k,m,l)$, if $P$ is true for the
considered value of $l$. Moreover,
if $P$ is false for the considered value of $l>2$, then
$b_l(G)=\infty$. In addition, one can require that the first two homology
groups of $G$ vanish.
\medskip\noindent
{\bf Proof of Proposition 3.3:}
 It is known how to effectively assign to any arithmetic
$\Sigma_3$ predicate $P$ of the form
$\exists n\forall m\exists k QQ(n,m,k)$
a Turing machine $t(P)$ of order $1$  such that
the halting set of $t(P)$ is cofinite if and only if
$P$ is true (cf. [So], p. 67). We are going to examine
this construction to
demonstrate that, in addition, the following assertion is true:
Assume that $P$ is true, and $N_0$ is the minimal value of $n$,
for which $\forall m\exists k\ QQ(n,m,k)$.
Then the halting set of $t(P)$ has a complement with
cardinality $N_0-1$.
A minor modification of this construction
will ensure that in the last case
the cardinality of the halting set will be not $N_0-1$ but
$N_0$.
\par
The construction of $t(P)$
consists of two steps. First, one replaces $\forall m\exists k
QQ(n,m,k)$ by a predicate asserting that a certain Turing machine $M(n)$
halts for infinitely many inputs: When $M(n)$ starts to work with input $m$
it checks all $p\leq m$, and for each $p$ it looks for
$k(p)$ that satisfies $QQ(n,m,k(p))$. The computation halts if and only if it
finds such $k$ for all $p\leq m$.  It is clear that the halting set
of $M(n)$ is infinite if and only if $M(n)$ halts with every input. And
this happens if
and only if for any $m$ there exists $k$ such that $QQ(n,m,k)$ holds. 
Note that if the halting set of $M(n)$ is infinite, then the set of values
of halting time is unbounded.
\par
The second step is slightly more complicated. We define
$t(P)$ by constructing its halting set $W$  or, more precisely, by
constructing the complement
of $W$.
The complement
to $W$ will be constructed in stages with the aid of infinitely many moving
markers numbered by $1,2, 3,\ldots$.
Think about numbers $1,2,\ldots$ as about being written on cells
of an infinite tape. Initially the
markers rest on all cells, so that the marker $i$ rests
on the cell number $i$. At the moment of time $s$ we check all markers
starting from the first in the increasing order of numeration until
the marker $s$. The marker
$i$ moves, if $i\leq s$, and
$s$ is the halting time of $M(i)$ with at least one of the
inputs $1,2,\ldots,s$. If the marker $i$ moves, then it moves to
the position occupied by the marker $(i+1)$, the marker $(i+1)$ moves to
the position occupied by the marker $(i+2)$, etc.
\par
Note that the movement
of markers $(i+1)$, $(i+2)$, etc. caused by the movement of the marker
$i$ is independent of their possible movement in the case
when $i+1$ (or $i+2$, etc.) 
turns out to be the halting time for $M(i+1)$ (correspondingly, $M(i+2)$)
and one of the inputs
$1,\ldots, s$.
\par
The cells that become free of
markers are then immediately enumerated to $W$. After infinitely
many steps markers will occupy all cells in the complement of $W$.
In other words, $W$ is the set of numbers of cells on the tape that will
be free of markers at some time.
It is clear that $W$ has a finite complement if and only if one of
the markers moves to infinity. If $N_0$ is the minimal number of a marker
that moves to infinity then the complement of $W$ will contain $N_0-1$
elements. (And a marker $i$ moves infinitely many times if and
only if $M(i)$ has an infinite halting set.) 
\par
Finally, to ensure that
in the last case the cardinality of ${\bf N}\setminus W$ is not $N_0-1$
but $N_0$,
we can use an infinite tape with cells numbered
by $0,1,2,\ldots$. The marker that stands at cell $0$ does not move, providing
us with a required extra element of the complement of $W$.
We add $1$ to every element that is being enumerated in $W$ in order
to return to ${\bf N}$ from ${\bf N}\bigcup \{0\}$.
QED.
\par\noindent
{\bf Remark 3.3.2.} A. In Proposition 3.3 and Corollary 3.3.1 we were assuming
that $n\in\{1,2,3\ldots\}$. Yet it is very easy to modify the proof of
Proposition 3.3 (and therefore of Corollary 3.3.1) for the case, when
$n\in\{0,1,2,\ldots\}$: One considers markers on an infinite tape cells of
which are numbered by $0,1,2\ldots$, but there 
is no ``dummy" (unmovable) marker at $0$.
\par\noindent
B. Note that the construction from [1] used in
the proof of Corollary 3.3.1 can be used to ensure
that the finitely presented group in the text of Corollary 3.3.1
has the following additional
property: Each of its homology groups is isomorphic to
the direct sum of a finite or infinite number of copies of ${\bf Z}$.
\medskip\noindent
{\bf Question 3.4.} Is it true that for
every partial function $f$
defined for $l=1, 2$ and computable by a Turing machine of
order $3$ there exists a finitely presented group
$G$ such that for every $l$ $b_l(G)=f(l)$, if $f(l)$ is defined, and $b_l(G)=\infty$, if $f(l)$ is not defined?
\medskip\noindent
Observe, that all these constraints on $f$ are necessary in order
for $f(n)$ to be the sequence of Betti numbers of a finitely
presented group. Therefore,
if the answer for this question is positive
one obtains a complete and very natural characterization of Betti numbers
of finitely presented groups.
A positive answer for Question 3.4 does not follow from Corollary 3.3.1
since we do not know how to effectively realize
$f(l)$ as the minimal $n$ such that
$\forall m\exists k QQ(n,m,k)$ is true. The construction that we used
ensures only that this minimal value of $n$ is greater than or equal 
to $f(l)$.
Observe, that by the
virtue of the discussion above the positive answer for this
question would follow from the positive answer for the
following question:
\medskip\noindent
{\bf Question 3.5.} Is there an algorithm that assigns
to every arithmetic predicate
$P=\exists n_1\exists n_2\forall m\exists k\
Q(n_1, n_2, m, k)$ a Turing machine $t$ such that
1) if $P$ is not true, then the halting set of $t$
has an infinite complement; 2) if $P$ is true, then the
complement of $t$ is finite and has cardinality equal to
the minimal value of $n_1$ such that $\exists n_2\forall m
\exists k\ Q(n_1,n_2,m,k)$.
\bigskip\noindent
{\bf 4. A finitely presented group with random Betti
numbers.}
\medskip\noindent
We recall the definition of Martin-L\"of randomness.
Let $2^\omega$ denote the set of all binary $0$-$1$ sequences identified
with $[0,1]$ interval via the binary representation. An {\it effective
null $G_\delta$ set} $S\subset 2^\omega$ is a countable intersection of a 
recursive sequence of $\Sigma_1^0$ subsets $U_n$, $n=1,2,\ldots$ of $S$
such that $\mu(U_n)\leq {1\over 2^n}$ for all $n$. (Here $\mu$ denotes
the Lebesque measure on $S=[0,1]$. A subset $U$ of $2^\omega$ is
$\Sigma_1^0$ if it can be represented as the set of all $f\in 2^\omega$ such
that $\exists n R(f,n)$, where $n$ runs over the set of
natural numbers, and $R$ is
a recursive predicate.) A sequence $f\in 2^\omega$ is {\it random}
if $f\not\in S$ for all effective null $G_\delta$ sets $S$.
\par
One can view effective null $G_\delta$ sets as effective randomness tests.
The intuitive idea behind this definition is that the sequence is random
when it passes all possible effective randomness tests.
For instance, $G=\{\{a_i\}\vert
\forall N \Sigma_{i=1}^N a_i<0.49 N\}$
is an effective null $G_\delta$ subset of
$2^\omega$, and because of the law
of large numbers a random binary sequence cannot be in this set.
\par
Another
(equivalent) definition of randomness is due to G. Chaitin. It uses
the notion of Kolmogorov complexity, i.e. the minimal length of a
description of objects by a program in a chosen model of computations.
Chaitin introduces the notion of self-delimiting Turing machines.
According to his definition a
binary sequence is random if and only if the Kolmogorov complexity
(=the minimal length of a description) of its first $n$ bits by a
self-delimiting Turing machine is $\geq n-c$ for some constant $c$. In
other words, a binary sequence is random if the sequences of
its first $n$ bits do not admit 
essentially better descriptions than just writing down all bits.
Martin-L\"of proved that the set of not random binary sequences
forms an effective null $G_\delta$ set (cf. Theorem 8.3 in [8]).
This result immediately implies that there exists a nonempty
$\Pi_1^0$ set $P\subset 2^\omega$ such that all its elements are random
(cf. Corollary 8.4 in [8]. A subset of $2^\omega$ is $\Pi_1^0$
if it can be represented as the
set of all $f\in 2^\omega$
such that $\forall n R(f,n)$, where $n$ runs over the set of natural
numbers, and $R$ is a recursive predicate.) This fact easily implies
that $P$ contains a binary sequence that can be presented as $\Delta_2$
predicate with three variables two of which are quantified and the third
is free. The meaning of the third variable,  $n$, is the number of a term
in the binary sequence: The $n$th term of the sequence is $1$ if and only if
the predicate is true for $n$. (Note that a predicate is
$\Delta_2$ if it is equivalent to a $\Pi_2$ predicate, and to a
$\Sigma_2$ predicate. We are grateful to Steve Simpson, who
explained to one of us this  sleek proof of the well-known fact
that there exist random $\Delta_2$ binary sequences.) We refer the reader
to [5], [2] for an introduction to random sequences and to [8] for a
modern treatment of this subject.
\par
It is known that the characteristic
function of a recursively enumerable set cannot be random,
so explicit descriptions of specific random sequences must be
quite indirect. Therefore, it is interesting that:
\medskip\noindent
{\bf Theorem 4.1.} There exists a finitely presented group
such that each of its
homology groups is either trivial or isomorphic to
{\bf Z}, and the sequence of Betti numbers of this group
is a random binary sequence.
\medskip\noindent
Indeed, it is known that there exists a subset of ${\bf N}$ 
defined by an arithmetic $\Delta_2$ predicate with one free
variable $n$ such that its
characteristic function represents a random binary sequence (see
the discussion before Theorem 4.1).
Now we
can write down a definition of this characteristic function 
as a $\Pi_2$ predicate with free variables $n$ and
$s$, where $s$ is set to be equal to one,
if the value of the characteristic function at $n$ is equal to one.
Denote this predicate by $R(s,n)$. Now consider the predicate
$P=\exists s R(s,n)$.
Note that $P$
is a $\Sigma_3$ predicate with a free variable $n$.
The minimal value of $s\in\{0,1,2\ldots\}$
for which $R(s,n)$ is true
is equal to the $n$th bit of the random binary sequence.
Now Theorem 4.1 follows from Corollary 3.3.1 and Remark 3.3.2.
QED.
\bigskip\noindent
{\bf 5. Even more rapidly growing functions arising in
homological group theory.}
\medskip\noindent
How rapidly can a function that naturally appears in ``mainstream" Mathematics (outside of Logic) can grow?
The following theorem demonstrates that
the fifth busy beaver functions can occur in a
natural way. On the other hand, one probably cannot
expect to find problems leading to appearance of
$B_k$ for significantly larger $k$ since such problems must
necessarily involve $k$ blocks of distinct quantifiers.
\medskip\noindent
{\bf Definitions:}
1. For every $N$ consider the set $C_N$ of all finite presentations of
groups $G$ of length $\leq N$ such that $b_i(G)$ is infinite
only for a finite set of indices $i$. For every group $G$ with a finite
presentation in $C_N$ denote the number of indices $j$ such that $b_j(G)$ is
infinite by $j(G)$. Define $c$ by the formula $c(N)=\max_{C_N}j(G)$.
\par\noindent
2. We say that a function $g$ from ${\bf N}$ to ${\bf N}$  {\it grows
as} a function $h:{\bf N}\longrightarrow {\bf N}$, if there exists
computable functions $\phi_1,\phi_2:{\bf N}\longrightarrow {\bf N}$
such that for all sufficiently large $n$ $g(\phi_1(n))>h(n)$ and
$h(\phi_2(n))>g(n)$. (Thus, all computable functions grow in the same
way, and according to Theorem 3.2 for every $k>2$ $b_k$ grows as $B_3$.
But $B_n$ does not grow as $B_m$, if $n\not= m$.)
\medskip\noindent
{\bf Theorem 5.1.} $c$ grows as $B_5$.
\medskip\noindent
{\bf Proof of Theorem 5.1:} The upper bound $c$
in terms of $B_5$ follows from the possibility to represent
sets $C_N$ in the definition of $c$ by
a $\Sigma_5$ predicate. To prove the 
upper bound for $B_5$ in terms of $c$
note the halting problem for Turing machines of
order $5$ can be effectively reduced to the problem $Cof$ for halting
sets of Turing machines of order $3$. (Recall that
$Cof$ is the algorithmic problem of
determination of cofiniteness.) The proof of the similar reducibility
for $Cof$ for recursively enumerable sets in [9] can be
easily generalized for halting sets of Turing machines of any order.
Moreover, following our analysis of this reducibility in the proof
of Theorem 3.1 we see that
there is an algorithm assigning to every Turing machine $T$ of order $5$
a Turing machine $t(T)$ of order $3$ such that $T$ halts with empty tape
if and only if the halting set of $t(T)$ has a finite complement, and if
the halting set of $t(T)$ has a finite complement, then its cardinality is
greater than or equal to the halting time of $T$.
\par
On the other hand the proof of Theorem 3.1
implies that there is an algorithm that for every Turing
machine of order $3$ finds a finitely presented group such
that for every $n\geq 1$ its Betti number $b_{n+2}$ is finite
if and only if the computation of this Turing machine of order $3$ with
input $n$ halts. Apply this construction to $t(T)$ and denote the resulting
finitely presented group by $G(T)$.
The number of indices $i$ for which $b_i(G(T))$ is
infinite will be equal to the cardinality of the complement
of the halting set of $t(T)$. The effectiveness
of the constructions implies that the length of the constructed
finite presentation of the group $G(T)$ is bounded by a computable
function of the number of states of $T$.
Taking the maximum over all Turing machines $T$ of order
$5$ that halt with the empty tape we obtain the desired
inequality. QED.
\bigskip\noindent
{\bf Acknowledgements.} This paper has been partially written during the
visit of 
Alexander Nabutovsky to
the Max Planck Institute for Mathematics in the Sciences in June-July, 2005.
Alexander Nabutovsky would like to thank the Max Planck Institute for
Mathematics in the Sciences for its kind hospitality.
\par
Both authors would like to acknowledge the partial support of this research
by NSF grants. A. Nabutovsky would like to acknowledge a partial support of this
research by NSERC Discovery grant.
\bigskip\noindent
{\bf References:}
\bigskip\noindent
[1] G. Baumslag, E. Dyer, C. Miller, ``On the integral homology of
finitely presented groups", Topology, 22(1983), 27-46.
\par\noindent
[2] G. Chaitin, ``Information and randomness", Springer-Verlag, 1994.
\par\noindent
[3] H. Enderton, ``Elements of recursion theory", in ``Handbook of
Mathematical Logic", ed. by J. Barwise, North-Holland, 1977, pp. 527-566.
\par\noindent
[4] D. Kan, W. Thurston, ``Every connected space has a homology of a $K(\pi, 1)$, Topology 15(1976), 253-258.
\par\noindent
[5] M. Li and P. Vitanyi, ``An introduction to Kolmogorov complexity and its
applications", Springer Verlag, 1997.
\par\noindent
[6] A. Nabutovsky, S. Weinberger, ``The fractal nature of Riem/Diff I",
Geom. Dedicata 101(2003), 1-54.
\par\noindent
[7] H. Rogers, Jr., ``Theory of recursive functions and effective computability", Second edition,
MIT Press, 1987.
\par\noindent
[8] S. Simpson, ``Mass problems and randomness",
Bull. Symbolic Logic, 11(1)(2005), 1-27.
\par\noindent
[9] R. Soare, ``Recursively-enumerable sets and degrees", Springer, 1987.
\par\noindent
[10] J. Stallings, ``A finitely presented group whose 3-dimensional integral
homology is not finitely generated", Amer. J. Math. 85(1963), 541-543. 
\bigskip\bigskip\noindent
\bye